\documentclass[12pt,a4paper]{article}

\usepackage{amsmath}
\usepackage{amssymb}
\usepackage{ifthen}

\usepackage{cases}
\newcommand{\ot}{\otimes}
\renewcommand{\a}{a_{ij}}
\newcommand{\G}{\Gamma}
\newcommand{\GH}{\hat{\G}}
\newcommand{\lam}{\lambda_{ij}}
\newcommand{\tth}{\tilde{\theta}}
\newcommand{\alp}{\alpha}
\newcommand{\si}{\sigma}
\newcommand{\ga}{\gamma}
\newcommand{\de}{\delta}
\newcommand{\bu}{$\bullet $ }

\newcommand{\pth}[1]{#1^{\text{th}}}

\newcommand{\dis}{\nsim}
\newcommand{\uf}{\mathfrak u}
\newcommand{\Uf}{\mathfrak U}
\newcommand{\D}{\mathcal D}
\newcommand{\X}{\mathcal X}
\newcommand{\A}{\mathcal A}
\newcommand{\U}{\mathcal U}
\newcommand{\B}{\mathcal B}
\newcommand{\K}{\mathcal K}
\newcommand{\Y}{\Upsilon}

\renewcommand{\k}{\Bbbk}

\newcommand{\liti}[6]{#1, \emph{#2}, #3 {\bf #4} (#5), #6. } 

\newtheorem{propo}{Proposition}
\newtheorem{thm}[propo]{Theorem}
\newtheorem{lemma}[propo]{Lemma}
\newcommand{\prf}[2][leer]{\begin{description} \ifthenelse{\equal{#1}{leer}}{\item[Proof:]}{\item[Proof #1:]} #2 \hspace{5em}  \hspace*{\fill}\bf{qed.} \end{description} } 

\DeclareMathOperator{\ad}{ad}
\DeclareMathOperator{\del}{\Delta}
\DeclareMathOperator{\eps}{\varepsilon}
\DeclareMathOperator{\an}{S}
\DeclareMathOperator{\ord}{ord}

\renewcommand{\_}[1]{_{\scriptscriptstyle (#1)}}

\usepackage[tips,graph,matrix]{xy}
\begin{document}

\author{Daniel Didt\footnote{This work will be part of the author's Ph.D. thesis written under the supervision of Professor H.-J. Schneider. The author is a member of the Graduiertenkolleg ``Mathematik im Bereich ihrer Wechselwirkung mit der Physik'' at Munich University.}\\\small Mathematisches Institut, LMU M\"unchen
, Theresienstr.~39\\\small 80333 M\"unchen, Germany\\\small didt@mathematik.uni-muenchen.de}
\title{Pointed Hopf algebras and Quasi-isomorphisms}
\date{}
\maketitle

\begin{abstract}
We show that a large class of finite dimensional pointed Hopf algebras is quasi-isomorphic to their associated graded version coming from the coradical filtration, i.e.~they are $2$-cocycle deformations of the latter. This supports a slightly specialized form of a conjecture in \cite{MA}.
\end{abstract}

\section{Introduction}
Recently there has been a lot of progress in determining the structure of pointed Hopf algebras. This has led to a discovery of whole new classes of such Hopf algebras and to some important classification results. See the survey article \cite{POIN} for an introduction and references.

A lot of these classes contain infinitely many non-isomorphic Hopf algebras of the same dimension, thus disproving an old conjecture of Kaplansky. Masuoka showed in \cite{MA}, and in a private note, that for certain of these new families the Hopf algebras are all $2$-cocycle deformations of each other. This led him to weaken Kaplansky's conjecture, stating that up to quasi-isomorphisms there should only be a finite number of Hopf algebras of a given dimension. This was disproved in \cite{EG} for a family of Hopf algebras of dimension 32. However, our results support Masuoka's conjecture in a really big class of examples. So we suggest that the conjecture could still be saved by specializing it slightly. We have the feeling that there is a fundamental difference between the even and odd dimensional case. We propose for a base field of characteristic zero:
\begin{quote} In a given odd dimension there are only finitely many non quasi-isomorphic pointed Hopf algebras whose coradical is commutative.\end{quote}
Two Hopf algebras are called \emph{quasi-isomorphic} or \emph{monoidally co-Morita equivalent} if their categories of comodules are monoidally equivalent. In \cite{S} Schauenburg showed that for finite dimensional Hopf algebras this is equivalent to the Hopf algebras being $2$-cocycle deformations of each other.
For the sake of completeness we quickly recall the notions associated with cocycle deformations.

Let $\A$ be a Hopf algebra over the base field $\k$. A linear, convolution invertible map $\si : \A\ot\A \rightarrow \k$ satisfying
\begin{align*}
\si(x\_1,y\_1)\si(x\_2y\_2,z)=\si(y\_1,z\_1)\si(x,y\_2z\_2)\quad
\text{and}\quad\si(1,1)&=1
\end{align*}
for all $x,y,z\in\A$ is called a \emph{unital 2-cocycle}. Given such a cocycle we can construct a new Hopf algebra $\A_{\si}$ which, as a coalgebra, is the same as $\A$ but has a new multiplication $\cdot_{\si}$ given by
\begin{equation*} x\cdot_{\si}y=\si(x\_1,y\_1)x\_2y\_2\si^{-1}(x\_3,y\_3).\end{equation*}
Here $\si^{-1}$ denotes the convolution inverse of $\si$. $A$ and $A_{\si}$ are then called \emph{cocycle deformations} of each other or \emph{quasi-isomorphic}.

We now want to describe the Hopf algebras we will be considering. We use the standard notation $\del,$ $\eps$ and $\an$ for the co-structures and leave out the summation sign in the Sweedler notation, i.e.~$\del(x)=x\_1\ot x\_2.$\\
We recall the notion of a \emph{linking datum $\D$ of finite Cartan type for a finite abelian group $\G$,} cf.~\cite{AS}. It consists of
\begin{itemize}
\item  a finite Cartan matrix $(a_{ij})_{1\le i,j\le\theta}$ of size $\theta\times\theta$.
\item a collection of group elements $g_1,\dots g_\theta\in\G$ and characters $\chi_1,\dots\chi_\theta\in\GH$ such that $N_i$, the order of $\chi_i(g_i)$, is odd, bigger than $3$ and not divisible by $3$ if $i$ belongs to a connected component of type $G_2$. Further we require \begin{equation}\label{cartype} \chi_i(g_j)\chi_j(g_i)=\chi_i(g_i)^{\a}\quad\text{for all}\;1\le i,j\le\theta. \end{equation}
\item a family $(\lam)_{1\le i,j\le\theta,i\dis j}$ of elements in the base field, where $\lam$ is arbitrary if $g_ig_j\neq 1$ and $\chi_i\chi_j=1,$ but $0$ otherwise. $i\dis j$ means vertices $i$ and $j$ lie in different connection components of the corresponding Dynkin diagram. Vertices with $\lam\neq 0$ are called \emph{linked}. In \cite[Lemma 5.6.]{AS} it is shown that a vertex can not be linked to two different vertices.\\ For technical reasons we require $\lam=-\chi_j(g_i)\lambda_{ji}$ and set $\lam=0$ when $i\sim j.$
\end{itemize}
For such a datum we can now consider the following algebra $\uf(\D)$. We fix a decomposition of the group $\G=<\!Y_1\!>\oplus\cdots\oplus <\!Y_s\!>$ and denote by $M_h$ the order of $Y_h$ for $1\le h\le s$. Then the algebra $\uf(\D)$ is generated by generators $a_1,\dots,a_\theta,y_1,\dots,y_s$ and relations
\begin{align}
 y_h^{M_h}=1,\quad y_hy_m&=y_my_h\qquad && \text{for all }1\le h,m\le s,\\
 ga_i&=\chi_i(g)a_ig\qquad && \text{for all }g\in\G,\,1\le i\le\theta,\\
 (\ad a_i)^{1-\a}a_j&=0\qquad && \text{for all }1\le i\neq j\le\theta, i\sim j\\
 \label{links}a_ia_j-\chi_j(g_i)a_ja_i&=\lam(1-g_ig_j)\qquad && \text{for all }1\le i\neq j\le\theta, i\dis j,\\
 \label{rootrel}a_\alp^{N_I}&=0\qquad && \text{for all }\alp\in\Phi^+_I, I\in\X.
\end{align}
Here $a_\alp$ denotes a so called root vector for a positive root $\alp$ of a connected component $I$ of the Dynkin diagram corresponding to the Cartan matrix $(\a)$. See also section $4.1.$ of \cite{AS}. $\X$ is the collection of all connected components and $N_I=N_i$ for any $i\in I$. This is well defined because of (\ref{cartype}).

$\ad a_i$ is the adjoint action defined for all $x\in\uf(\D)$ by $$(\ad a_i)x:=a_ix-g_ixg_i^{-1}a_i.$$
Elements $g\in\G$ are naturally interpreted as words in the generators $y_i$.

Theorem 5.17.~of \cite{AS} now states that $\uf(\D)$ is a pointed Hopf algebra of dimension $|\G|\prod_{I\in\X}N_I^{|\Phi^+_I|}$ with coradical equal to $\k\G.$
The Hopf algebra structure is determined by
\begin{equation*} \del(y_h)=y_h\ot y_h,\qquad \del(a_i)=a_i\ot 1+g_i\ot a_i.\end{equation*}
The associated graded Hopf algebra of $\uf(\D)$ coming from the coradical filtration is simply $\uf(\D_0),$ where $\D_0$ is the same as $\D$ except that all $\lam$ are zero.\\

Actually, one can even extend this class of Hopf algebras by replacing the zero in (\ref{rootrel}) with appropriately chosen central elements $u_{\alp}$ from the group algebra $\k\G.$ This however, has been done explicitly only in the cases where the given Cartan matrix is of type $A_n$ \cite{POIN}, $B_2$ \cite{BDR} or copies of $A_1$ \cite{P^3}. In each of these cases there is an infinite family of such possible central elements.

So we see that after fixing the abelian group $\G,$ certain elements $g_i,\chi_i$ and the Cartan matrix $\a$ we get families of Hopf algebras of the same dimension that have two sets of parameters: the linking parameters $\lam$ and the root vector parameters $u_{\alp}.$ All possible choices for these parameters can be determined. For the linking parameters, all the necessary considerations are in \cite{D}. In the cases where the root vector parameters are known explicitly, all of them have been found.

Classification results for classes of finite dimensional pointed Hopf algebras with abelian coradical of odd dimension seem to indicate that apart from a few exceptional cases all such Hopf algebras are of the above form.\\

Our ultimate goal is to show, that in each such family there is only one quasi-isomorphism class of Hopf algebras. This means that all Hopf algebras which differ only in their choice of linking and root vector parameters are cocycle deformations of each other. In this article we present three major steps towards our goal.

First we prove the statement in the case with only linking parameters and no root vector parameters, i.e.~(\ref{rootrel}) holds.

In the second part we show the result for the case where the Cartan matrix is of type $A_n.$ Here all the root vector parameters are known, and linking parameters do not appear.

In the last part we combine all results to treat the mixed case, where the Dynkin diagram is a union of $A_n$'s.

Once the root vector parameters for the other Cartan matrices of finite type have been determined explicitly, an analogous treatment should provide the same result as for the $A_n$ case.

\section{The linking case}

\begin{thm}\label{linkthm} For two linking data $\D$ and $\D'$ of finite Cartan type for a finite abelian group $\G$ which only differ in their choice of the $\lam$ we have that $\uf(\D)$ and $\uf(\D')$ are quasi-isomorphic.
\end{thm}

\prf{The strategy of the proof is as follows. We will prove the statement for an arbitrary linking datum and the datum where all the $\lam$ are 0.
This will be achieved inductively by proving the statement for an arbitrary linking datum with at least one pair of linked vertices and one datum with the same data, except that the number of connected components that are not linked to any other vertex is increased by one. This means, if in the original datum we have $\lam\neq 0$ for some $i,j$ then we take for the other datum $\lambda_{kl}=0$ for all $1\le l\le\theta$ and for all $k$ that are in the same connection component $I$ as $i$. The transitivity of the quasi-isomorphism relation will then yield the result.

So let $\D=\{(\a)_{1\le i,j\le\theta},(g_i)_{1\le i\le\theta},(\chi_j)_{1\le j\le\theta},(\lam)_{1\le i<j\le\theta}\} $ be a linking datum, where there is a connected component $I$ and an $i\in I$ such that there is a $j$ with $\lam\neq 0$. Let $\tilde{\D}$ be the same linking datum except that $\lambda_{kl}=0$ for all $k\in I$. It suffices to show that $\uf(\D)$ and $\uf(\tilde{\D})$ are quasi-isomorphic.

Now, from the proof of Theorem 5.17.~in \cite{AS} we see that $\uf(\D)$ is isomorphic to $(\U\ot\B)_{\si}/\K^+,$ where $\U,\B$ and $\K$ are Hopf algebras to be described shortly. $\si$ is a 2-cocycle constructed from the linking datum. The key observation is now that $\uf(\tilde{\D})$ is isomorphic to $(\U\ot\B)_{\tilde{\si}}/\K^+$ with the \emph{same} Hopf algebras involved and only the cocycle $\tilde{\si}$ is different. 

We reorder the vertices so that $I=\{1,\dots,\tilde{\theta}\}$ and set $\Y:={<\!Z_1\!>}\oplus\cdots\oplus {<\!Z_{\tth}\!>},$ where the order of $Z_i$ is the least common multiple of $\ord g_i$ and $\ord \chi_i.$ Let $\eta_j$ be the unique character of $\Y$ such that $\eta_j(Z_i)=\chi_j(g_i),$ $1\le i,j\le\tilde{\theta}.$ This is well defined because $\ord g_i$ divides $\ord Z_i$ for all $i.$

Now take as linking datum $\D_1$ for the original group $\G$ $$\D_1=\{(\a)_{\tilde{\theta}<i,j\le\theta},(g_i)_{\tilde{\theta}<i\le\theta},(\chi_j)_{\tilde{\theta}<j\le\theta},(\lam)_{\tilde{\theta}<i<j\le\theta}\}$$ and as $\D_2$ for the group $\Y$
$$\D_2=\{(\a)_{1\le i,j\le\tilde{\theta}},(Z_i)_{1\le i\le\tilde{\theta}},(\eta_j)_{1\le j\le\tilde{\theta}},(\lam)_{1\le i<j\le\tilde{\theta}}\}.$$
Then $\B:=\uf(\D_1)$ and $\U:=\uf(\D_2).$ We denote the generators of $\B$ by $b_{\tth+1},\dots,b_{\theta}$ and $y_1,\dots,y_s$ and the generators of $\U$ by $u_1,\dots,u_{\tth}$ and $z_1,\dots,z_{\tth}.$

The central Hopf subalgebra $\K$ of $\U\ot\B$ is $\k[z_i\ot g_i^{-1}:1\le i\le\tth].$

The cocycle $\si$ for $\U\ot\B$ is defined by \begin{equation*} \si(u\ot a,v\ot b):=\eps_{\U}(u)\tau(v,a)\eps_{\B}(b).\end{equation*}
$\tau: \U\ot\B\rightarrow\k$ is a linear map with the following list of properties
\begin{enumerate}
\item $\tau(uv,a)=\tau(u,a\_1)\tau(v,a\_2)$
\item $\tau(u,ab)=\tau(u\_1,b)\tau(u\_2,a)$
\item $\tau(1,a)=\eps_{\B}(a)$
\item $\tau(u,1)=\eps_{\U}(u).$
\end{enumerate}
It is given by $$\tau(u\ot b):=\varphi(u)(b),$$ where $\varphi:\U\rightarrow(\B^*)^{cop}$ is a Hopf algebra homomorphism defined on the generators of $\U$ by
$$\varphi(z_i):=\ga_i\qquad\text{and}\qquad\varphi(u_j):=\de_j.$$
Here $\ga_i:\B\rightarrow\k$ is a character defined on the generators of $\B$ by
$$ \ga_i(y_k):=\chi_i(y_k)\qquad\text{and}\qquad\ga_i(b_l):=0,$$
and $\de_j:\B\rightarrow\k$ is a $(\eps_{\B},\ga_j)$-derivation defined by
$$ \de_j(y_k):=0\qquad\text{and}\qquad\de_j(b_l):=-\chi_j(g_l)\lambda_{jl}.$$
In all the above formulas $1\le i,j\le\tth$ and $\tth<k,l\le\theta.$ The cocycle $\tilde{\si}$ is now defined in exactly the same way as $\si$ only in the last part $\tilde{\de}_j(b_l):=0$ as all the $\lambda_{jl}=0$ because we wanted for $\tilde{\D}$ the component $I=\{1,\dots,\tth\}$ not to be linked to any other component.

The inverse $\si^{-1}$ is given in the same way by $\tau^{-1},$ where $$\tau^{-1}(u,b):=\varphi(\an u)(b)=\varphi(u)(\an^{-1}b).$$
If $\si,\tilde{\si}$ are two 2-cocycles for a Hopf algebra $\A,$ then $\rho:=\tilde{\si}\si^{-1}$ is again a 2-cocycle, but for the Hopf algebra $\A_{\si}.$ In our case this means $\rho$ is a 2-cocycle for the Hopf algebra $\A:=(\U\ot\B)_{\si}.$ Then $(\U\ot\B)_{\tilde{\si}}=\A_{\rho}$. If we can show that $\rho$ passes down naturally to a 2-cocycle $\rho'$ on $\A/\K^+$ such that $(\A/\K^+)_{\rho'}$ is isomorphic to $\A_{\rho}/\K^+$ then the statement is clear:
\begin{equation} \uf(\tilde{\D})\simeq \A_{\rho}/\K^+\simeq (\A/\K^+)_{\rho'}\simeq \uf(\D)_{\rho'}. \end{equation}
Hence we want the following situation.
$$
\xy \xymatrix@C+1cm@R+1cm{\A\ot\A  \ar[r]^{\rho} \ar[d]^{\pi} & \k \\
     (\A\ot\A)/(\A\ot\K^++\K^+\ot\A) \simeq \A/\K^+\ot\A/\K^+ \ar@{{}-->}[ru]^{\rho'}\\} \endxy
$$
For this it suffices to show that $\rho$ is 0 on the kernel of the natural projection $\pi$ and we get the factorization and $(\A/\K^+)_{\rho'}\simeq \A_{\rho}/\K^+$ by definition of $\rho'.$
So we see that it is enough to show $\rho(\A,\K^+)=0=\rho(\K^+,\A).$ This means that for all $u\in\U,$ $b\in\B$ and $1\le i\le\tth$ we need
\begin{align*} \rho(u\ot b,z_i\ot g_i^{-1})&=\rho(u\ot b,1\ot 1)\\
                \rho(z_i\ot g_i^{-1},u\ot b)&=\rho(1\ot 1,u\ot b).
\end{align*}
We calculate using the definition of $\rho$, the convolution product and property 3.~of $\tau$ and $\tilde{\tau}:$
\begin{align*}\rho(u\ot b,z_i\ot g_i^{-1})&=\tilde{\si}(u\_1\ot b\_1,z_i\ot g_i^{-1})\si^{-1}(u\_2\ot b\_2,z_i\ot g_i^{-1})\\
                &=\eps_{\U}(u\_1)\tilde{\tau}(z_i,b\_1)\eps_{\B}(g_i^{-1})\eps_{\U}(u\_2)\tau(\an z_i,b\_2)\eps_{\B}(g_i^{-1})\\
                &=\eps_{\U}(u)\tilde{\tau}(z_i,b\_1)\tau(z_i^{-1},b\_2)\\
        \rho(u\ot b,1\ot 1)&=\tilde{\si}(u\_1\ot b\_1,1\ot 1)\si^{-1}(u\_2\ot b\_2,1\ot 1)\\
                &=\eps_{\U}(u\_1)\tilde{\tau}(1,b\_1)\eps_{\B}(1)\eps_{\U}(u\_2)\tau(\an 1,b\_2)\eps_{\B}(1)\\
                &=\eps_{\U}(u)\eps_{\B}(b).
\end{align*}
As $z_i$ is group-like we see by using property 2.~of $\tau$ and $\tilde{\tau}$ that it is enough to verify 
\begin{equation} \label{chk1}\tilde{\tau}(z_i,b\_1)\tau(z_i^{-1},b\_2)=\eps_{\B}(b) \end{equation} on the generators of $\B.$
An analog calculation, using property 1.~and 4.~this time, shows that for the second condition we have to verify for the generators of $\U$
\begin{equation} \label{chk2}\tilde{\tau}(u\_1,g_i^{-1})\tau(u\_2,g_i)=\eps_{\U}(u). \end{equation}
The verification goes as follows ($1\le j\le\tth, \tth<k\le\theta$)
\begin{align*}
\intertext{$b=y_k$ in (\ref{chk1})}
  \tilde{\tau}(z_i,y_k)\tau(z_i^{-1},y_k)&=\tilde{\ga}_i(y_k)\ga_i^{-1}(y_k)\\
                                &=\chi_i(y_k)\chi_i^{-1}(y_k)=1=\eps_{\B}(y_k),\\
\intertext{$b=b_k$ in (\ref{chk1})}
  \tilde{\tau}(z_i,b_k)\tau(z_i^{-1},1)+\tilde{\tau}(z_i,g_k)\tau(z_i^{-1},b_k)&=\tilde{\ga}_i(b_k)\ga_i^{-1}(1)+\tilde{\ga}_i(g_k)\ga_i^{-1}(b_k)\\
                                        &=0+0=\eps_{\B}(b_k),\\
\intertext{$u=z_j$ in (\ref{chk2})}
 \tilde{\tau}(z_j,g_i^{-1})\tau(z_j,g_i)&=\chi_j(g_i^{-1})\chi_j(g_i)=\chi_j(1)\\
                                &=1=\eps_{\U}(z_j),\\
\intertext{$u=u_j$ in (\ref{chk2})}
 \tilde{\tau}(u_j,g_i^{-1})\tau(1,g_i)+\tilde{\tau}(z_j,g_i^{-1})\tau(u_j,g_i)&=\tilde{\de}_j(g_i^{-1})\eps_{\B}(g_i)+\tilde{\ga}_j(g_i^{-1})\de_j(g_i)\\
                                &=0+0=\eps_{\U}(u_j).
\end{align*}
}
\section{A special Root vector case}

Now we deal with Hopf algebras where the Dynkin diagram is just one copy of $A_n,$ $n\ge 1,$ and hence there is no linking. To fix the order of the indices we require for the $n\times n$ Cartan matrix $a_{ij}=0$ whenever $|i-j|\ge 2.$
Given an abelian group $\G$ and a corresponding linking Datum $\D$ of $A_n$ type, we consider again the Hopf algebra $\uf(\D)$ only we replace $0$ in the root vector relation (\ref{rootrel}) by certain central elements $u_{\alp}$ from the group algebra.\\
To be more precise, we consider now the following setup.\\

$R$ is a left $\k\G$-module algebra generated by $a_1,\dots,a_n,$ subject to the quantum Serre relations
\begin{equation} (\ad a_i)^{1-\a}a_j=0. \end{equation}
The action of $g\in\G$ is given by $g.a_i=\chi_i(g)a_i$ and the adjoint action can be expressed as $(\ad a_i)x=a_ix-(g_i.x)a_i.$
We denote by $H_0:=R\#\k\G$ the smash product. It is a Hopf algebra with the comultiplication given by $\del(a_i)=a_i\ot 1+g_i\ot a_i.$

In \cite{POIN} the explicit form of the root vectors $a_{\alp}$ and all the possible families $u_{\alp}$ have been determined. We repeat these here for later calculations.

Set $q_{i,j}:=\chi_j(g_i),\;1\le i,j\le n,$ $N:=\ord \chi_i(g_i)$ and $B^{i,j}_{p,r}:=\prod_{i\le l<j, p\le h<r}q_{l,h}.$ $N$ is defined independently of the choice of $i$ as $\chi_i(g_i)=\chi_j(g_j)$ due to (\ref{cartype}). Then $$C^{j}_{i,p}:=(1-q^{-1})^N(B^{p,j}_{i,p})^{\binom N2},\qquad 1\le i<p<j\le n+1.$$
The root vectors are defined inductively on the height of the corresponding root:
\begin{align}
e_{i,i+1}&:=a_i\;,\qquad 1\le i\le n;\\
\label{defroot}e_{i,j}&:=e_{i,j-1}e_{j-1,j}-B^{i,j-1}_{j-1,j}e_{j-1,j}e_{i,j-1}\qquad j-i\ge 2.
\end{align}
We set $$ \chi_{i,j}:=\prod_{i\le l<j}\chi_l,\quad g_{i,j}:=\prod_{i\le l<j}g_l\quad\text{and}\quad h_{i,j}:=g_{i,j}^N.$$
Then $B^{i,j}_{p,r}=\chi_{p,r}(g_{i,j}).$
Next we introduce a parameter family $\ga$ containing for every root vector $e_{i,j}$ a scalar $\ga_{i,j}\in\k.$ The elements $u_{i,j}\in\k\G$ are also defined inductively:
\begin{equation}\label{defu} u_{i,j}(\ga):=\ga_{i,j}(1-h_{i,j})+\sum_{i<p<j}C^j_{i,p}\ga_{i,p}u_{p,j}(\ga),\quad 1\le i<j\le n+1. \end{equation}
We require that for all $i<j,$
\begin{equation}\label{charn}\ga_{i,j}=0\quad\text{if}\quad\chi_{i,j}^N\neq\eps\;\text{or}\; h_{i,j}=1.\end{equation}
We will call such families \emph{admissible}. In \cite[Lemma 7.20.]{POIN} it was shown that this condition implies the same statement for the $u_{i,j}$ and proves all the $u_{i,j}$ to be central in $H_0.$ 

Finally, we define the Hopf algebra $\A(\ga)$ as the quotient of $H_0$ by the ideal generated by the root vector relations
\begin{equation} e_{i,j}^N=u_{i,j}(\ga),\qquad 1\le i<j\le n+1. \end{equation}
Our result is now:
\begin{thm}\label{main}For a linking datum of $A_n$ type and two admissible families $\ga$ and $\ga'$ of scalars the above described Hopf algebras $\A(\ga)$ and $\A(\ga')$ are quasi-isomorphic.
\end{thm}
For the proof of this theorem we will first slightly generalize Theorem 7.24.~in \cite{POIN}.

Let $\G$ be an abelian group and $R$ a left $\k\G$-module algebra, $H_0:=R\#\k\G.$ Assume that there is an integer $P$, a subset $Z$ of $\{1,2,\dots,P\},$ integers $N_j>1$ where $j\in Z$ and elements $y_i\in R,$ $h_i\in\G,$ $\eta_i\in\GH,$ $1\le i\le P,$ such that
\begin{align*}
g\cdot y_i&=\eta_i(g)y_i\quad\text{for all }g\in\G,\; 1\le i\le P.\\
y_iy_j^{N_j}&=\eta_j^{N_j}(h_i)y_j^{N_j}y_i\quad\text{for all }1\le i\le P,\;j\in Z.\\
\text{The ele}&\text{ments }y_1^{a_1}\dots y_P^{a_P},\;a_1,\dots a_P\ge 0,\text{ form a }\k\text{-basis of }R.
\end{align*}
\begin{thm}\label{basisthm}
Let $u_j,\;j\in Z,$ be a family of elements in $\k\G,$ and $I$ the ideal in $H_0$ generated by all $y_j^{N_j}-u_j,$ $j\in Z.$ Let $A=H_0/I$ be the quotient algebra.

If $u_j$ is central in $H_0$ for all $j\in Z$ and $u_j=0$ if $\eta_j^{N_j}\neq\eps$,\\
then the residue classes of $y_1^{a_1}\dots y_P^{a_P}g,$ $a_i\ge 0,$ $a_j<N_j$ if $j\in Z,$
$g\in\G,$ form a $\k$-basis of $A.$
\end{thm}
The proof of this theorem is exactly the same as in the original paper, where $Z$ included all indices from $1$ to $P.$\\

To see how this can be applied in our situation, we recall two more results (Theorem 7.21.~and Lemma 7.22.) from \cite{POIN}.
First, the elements $e_{1,2}^{b_{12}}e_{1,3}^{b_{13}}\dots e_{n,n+1}^{b_{nn+1}}g,$ $g\in\G,$ $b_{ij}\ge 0,$ where the root vectors are arranged in the lexicographic order, form a $\k$-basis of $H_0.$ Furthermore, we have the following crucial commutation rule for all $1\le i<j\le n+1,$ $1\le s<t\le N+1$
\begin{equation}\label{crucial}
e_{i,j}e_{s,t}^N=\chi_{s,t}^N(g_{i,j})e_{s,t}^Ne_{i,j}.
\end{equation}
So the last theorem gives us, for instance, a basis of $A(\ga).$

We are now ready to give the proof of our Theorem.
\prf[of Theorem \ref{main}]{As in the linking case, it is enough to prove that for any admissible family $\ga$, $A(\ga)$ is quasi-isomorphic to $A(\ga_0)$ where $\ga_0$ denotes the family where all the $\ga_{i,j}$ are zero. Again, this  will be achieved by following a stepwise procedure.

Fix $i_0$ with $1\le i_0\le n$ such that for the given admissible family $\ga$ we have $$\ga_{i,j}=0 \text{ for all }1\le i<i_0\text{ and }i<j\le n+1.$$ Set $\tilde{\ga}_{i,j}:=\ga_{i,j}$ if $i\neq i_0$ and $\tilde{\ga}_{i_0,j}:=0$ for all $i_0<j\le n+1.$ Then $\tilde{\ga}$ is again admissible. It is now sufficient to prove that $\A(\ga)$ and $\A(\tilde{\ga})$ are quasi-isomorphic and then to repeat this step with increasing $i_0,$ replacing $\ga$ by $\tilde{\ga}.$
Set
$$H=H_0/I,\; I:=(e_{i,j}^N-u_{i,j}\,:\,1\le i<j\le n+1,\;i\neq i_0\text{ or }\chi_{i,j}^N\neq\eps).$$
Note that $u_{i,j}(\ga)=u_{i,j}(\tilde{\ga})$ for all the $u_{i,j}$ appearing in the ideal $I.$

$H$ is a Hopf algebra.\\
To see this, we recall the comultiplication on the root vectors.
\begin{equation}\begin{split}\label{delroot}&\del(e_{i,j}^N-u_{i,j})=(e_{i,j}^N-u_{i,j})\ot 1+h_{i,j}\ot(e_{i,j}^N-u_{i,j})+\\
&+\sum_{i<p<j}C_{i,p}^je_{i,p}^Nh_{p,j}\ot(e_{p,j}^N-u_{p,j})+\sum_{i<p<j}C_{i,p}^j(e_{i,p}^N-u_{i,p})h_{p,j}\ot u_{p,j}
\raisetag{4em}
\end{split}\end{equation}
If $i<i_0$ we have $\ga_{i,j}=0$ for all $i<j\le n+1$ and hence $u_{i,j}=0.$ So all the summands in (\ref{delroot}) are in $H_0\ot I$ or $I\ot H_0,$ because $e_{i,p}^N\in I$ for all $i<p<j.$\\
The case $i>i_0$ is obvious. When $i=i_0$ then, according to the definition of $I$, $e_{i_0,j}^N-u_{i_0,j}$ is in $I$ only if $\chi_{i_0,j}^N\neq\eps.$ From (\ref{charn}) follows then $\ga_{i_0,j}=0$ and $u_{i_0,j}=0.$ We have 
$$\chi_{i_0,j}^N=\chi_{i_0,p}^N\chi_{p,j}^N\qquad\text{for all }i_0<p<j.$$
And hence $\chi_{i_0,p}^N\neq\eps$ and $u_{i_0,p}=0$ or $\chi_{p,j}^N\neq\eps$ and $u_{p,j}=0.$ This proves that 
$$\del(I)\subset H_0\ot I+I\ot H_0.$$
A simple calculation shows $\eps(e_{i,j})=0=\eps(u_{i,j}).$ So we can get two recursion formulas for the antipode of $e_{i,j}^N-u_{i,j}$ from the comultiplication formula, as $\an(a\_1)a\_2=\eps(a)=a\_1\an(a\_2).$ Using the second of these formulas for the case $i<i_0$ we have immediately $\an(e_{i,j}^N)\in I.$ When $i>i_0$ an inductive argument using the first formula and $$\an(e_{i,i+1}^N-u_{i,i+1})=-g_i^{-N}(e_{i,i+1}^N-u_{i,i+1})$$ gives again $\an(e_{i,j}^N)\in I.$ For $i=i_0,$ a combination of the reasoning from the discussion of the comultiplication and the inductive argument from the last case give the desired result, establishing $\an(I)\subset I.$ Hence $I$ is a Hopf ideal.

Next we take as $K$ the sub algebra of $H$ generated by the group $\G$ and the remaining $e_{i_0,j}^N,$ i.e.~$i_0<j\le n+1,$ such that $\chi_{i_0,j}^N=\eps.$ A similar calculation as above reveals that $K$ is actually a sub Hopf algebra of $H.$ This time one has to use the fact that for any $p$ between $i_0$ and $j$ either $\chi_{i_0,p}^N=\eps$ or that $\chi_{p,j}^N\neq\eps$ and hence $u_{p,j}=0=e_{p,j}^N$ in $H.$

As all the $u_{i,j}$ fulfill the conditions of Theorem \ref{basisthm} we immediately get a basis of $H.$ The commutation relations (\ref{crucial}) for the $\pth N$ powers of the root vectors show that the generators of $K$ all commute with each other, because the factor is $1,$ as $\chi_{i_0,j}^N=\eps$ for all generators. Hence any monomial in $K$ can be reordered and is then a basis element of $H.$ So $K$ is just the polynomial algebra on its generators.

We define an algebra map $f:K\rightarrow \k$ by setting 
$$f(e_{i_0,j}^N):=\ga_{i_0,j},\qquad f(g):=1\qquad\text{on all the generators of }K,\,g\in\G.$$
Algebra maps from a Hopf algebra $K$ to the base field form a group under the convolution product, where the inverse is given by the composition with the antipode. This group acts on the Hopf algebra $K$ from the left and the right by
$$f.x=x\_1f(x\_2),\qquad x.f=f(x\_1)x\_2.$$
To be able to apply Theorem 2.~of \cite{MA}, which will give us the desired quasi-isomorphism, we have to calculate $f.e_{i_0,j}^N.f^{-1}.$ As a preparation for this we calculate first $f(u_{i,i+1})=0$ for all $1\le i\le n$ and see then from the inductive definition (\ref{defu}) that $f(u_{i,j})=0$ for all $i<j.$ So for the generators of $K$ we have
\begin{equation*}\begin{split} f.e_{i_0,j}^N&=e_{i_0,j}^Nf(1)+h_{i_0,j}f(e_{i_0,j}^N)+\sum_{i_0<p<j}C_{i_0,p}^je_{i_0,p}^Nh_{p,j}f(e_{p,j}^N)\\&=e_{i_0,j}^N+\ga_{i_0,j}h_{i_0,j},
\end{split}\end{equation*}
as $e_{p,j}^N=u_{p,j}$ in $H.$ The recursive formula for the antipode of the $\pth N$ powers of the root vectors is 
\begin{equation*} \an(e_{i_0,j}^N)=-h_{i_0,j}^{-1}e_{i_0,j}^N-\sum_{i_0<p<j}C_{i_0,p}^jh_{p,j}^{-1}\an(e_{i_0,p}^N)e_{p,j}^N.
\end{equation*}
So we get
\begin{equation}\begin{split}\label{fif} (f.&e_{i_0,j}^N).f^{-1}=\ga_{i_0,j}f(h_{i_0,j}^{-1})h_{i_0,j}+f(\an(e_{i_0,j}^N))+f(h_{i_0,j}^{-1})e_{i_0,j}^N+\\ 
&\qquad\qquad\qquad+\sum_{i<p<j}C_{i,p}^jf(h_{i_0,j}^{-1})f(\an(e_{i_0,p}^N))e_{p,j}^N\\
&=e_{i_0,j}^N+\ga_{i_0,j}h_{i_0,j}+\\
&\quad+f\left(-h_{i_0,j}^{-1}e_{i_0,j}^N-\!\sum_{i_0<p<j}C_{i_0,p}^jh_{p,j}^{-1}\an(e_{i_0,p}^N)e_{p,j}^N\right)+\\
&\quad+\sum_{i<p<j}C_{i,p}^jf\left(-h_{i_0,p}^{-1}e_{i_0,p}^N-\!\sum_{i_0<q<p}C_{i_0,q}^ph_{q,p}^{-1}\an(e_{i_0,q}^N)e_{q,p}^N\right)u_{p,j}\\
&=e_{i_0,j}^N+\ga_{i_0,j}h_{i_0,j}-f(e_{i_0,j}^N)-\!\sum_{i<p<j}C_{i,p}^jf(e_{i_0,p}^N)u_{p,j}\\
&=e_{i_0,j}^N-\ga_{i_0,j}(1-h_{i_0,j})-\!\sum_{i<p<j}C_{i,p}^j\ga_{i_0,p}u_{p,j}\\
&=e_{i_0,j}^N-u_{i_0,j}\qquad\text{by (\ref{defu})}.
\end{split}\end{equation}
Two parts of the sum vanish as they are multiples of $f(e_{p,j}^N=u_{p,j})=0=f(e_{q,p}^N=u_{q,p}).$ Note that in the second last step, if $j>i_0$ is such that $\chi_{i_0,j}^N\neq\eps,$ then $e_{i_0,j}^N$ is not a generator of $K$ and we can not just apply the definition of $f.$ But in this case $u_{i_0,j}=0=e_{i_0,j}^N$ in H and $\ga_{i_0,j}$ is zero, too, as this is required for an admissible family. So we still have $f(e_{i_0,j}^N)=\ga_{i_0,j}.$

Let $J$ be the Hopf ideal of $K$ generated by all the generators $e_{i_0,j}^N$ of $K.$ Then,
according to \cite[Theorem 2.]{MA}, $H/(f.J)$ is an $(H/(J),H/(f.J.f^{-1}))$-bi-Galois object and hence $H/(J)$ and $H/(f.J.f^{-1})$ are quasi-isomorphic if the bi-Galois object is not zero.
We see that $u_{i_0,j}(\tilde{\ga})=0$ and so $A(\tilde{\ga})=H/(J).$ Calculation (\ref{fif}) showed that $A(\ga)=H/(f.J.f^{-1}).$ We are left to show that $B:=H/(f.J)$ is not zero.

$B=H_0/(I,f.J)$ by construction. We have a basis of $H_0$ and see that we could apply Theorem \ref{basisthm} to get a basis of $B$. It just remains to check that the elements $\ga_{i_0,j}h_{i_0,j}$ appearing in $f.J$ satisfy the conditions of the theorem.

If $\chi_{i_0,j}^N\neq\eps$ then $\ga_{i_0,j}=0$, because $\ga$ is admissible. $h_{i_0,j}$ is in $\G$ and so commutes with all group elements. We will show now that $h_{i_0,j}$ commutes also with all the generators of $H_0.$ For this we calculate
\begin{equation*}\begin{split}h_{i_0,j}a_k&=\chi_k(h_{i_0,j})a_kh_{i_0,j},\\
\chi_k(h_{i_0,j})&=\chi_k(\prod_{i_0\le p<j}g_p^N)
=\prod_{i_0\le p<j}\chi_k^N(g_p)\\&=\prod_{i_0\le p<j}\chi_p^N(g_k^{-1})
=\chi_{i_0,j}^N(g_k^{-1})=1.
\end{split}\end{equation*}
Here we used (\ref{cartype}) and the fact that $N$ is exactly the order of any diagonal element $\chi_p(g_p).$ So $B$ is not zero and the statement is proven.}

\section{The mixed case}

Now finally, we consider the case where linking parameters \emph{and} root vector parameters appear. As we do not know yet, how to generalize our considerations to arbitrary Dynkin diagrams of finite type, we will consider only copies of $A_n.$

So we take a linking datum $\D$ of type $A_{n_1}\times\dots\times A_{n_t},$ $n_k,t>0,$ for a fixed finite abelian group $\G.$
The algebra $\Uf(\D)$ is defined in the same way as $\uf(\D)$ in the introduction except that we leave out relation (\ref{rootrel}). This is a Hopf algebra, cf. \cite{D,New-AS}. We order the vertices in the Dynkin diagram so that the root vectors of the $\pth k$ component are $e_{S_k+i,S_k+j},$ $1\le i<j\le n_k+1,$ 
where $S_k=n_1+\dots+n_{k-1}.$ The root vectors within one component are defined in the same way as in the previous section. Then the monomials
$$ e_{1,2}^{b_{1,2}}\cdots e_{1,n_1+1}^{b_{1,n_1+1}}e_{2,3}^{b_{2,3}}\cdots e_{S_t+n_t,S_t+n_t+1}^{b_{S_t+n_t,S_t+n_t+1}},\quad 0\le b_{S_k+i,S_k+j}\le N_k,$$
form a PBW-basis of $\Uf(\D).$ Here $N_k$ is again the common order of the diagonal elements $\chi_i(g_i)$ with $i$ in the $\pth k$ component of the diagram. A proof can be found in \cite[Theorem 4.2.]{New-AS}. The idea is that a PBW-basis is known for the linking datum $\D_0$ where all the $\lam$ are zero. Using, for instance, the considerations in the second section of this paper one knows explicitly the 2-cocycle relating $\Uf(\D)$ and $\Uf(\D_0).$ Expressing now the above monomials in $\Uf(\D_0)$ one sees that they are basis elements of the common underlying vector space.

As in the third section we now introduce for every component $k$ of the diagram an admissible parameter family $\ga_k$ and define the corresponding elements $u_{S_k+i,S_k+j}(\ga_k)\in\k\G,$ $1\le i<j\le n_k+1.$ The collection of all the parameters $\ga_k$ will be denoted by $\ga.$

$\A(\D,\ga)$ is now defined as the quotient of $\Uf(\D)$ by the ideal generated by the root vector relations
$$ e_{S_k+i,S_k+j}^{N_k}=u_{S_k+i,S_k+j}(\ga_k),\qquad 1\le k\le t,\,1\le i<j\le n_k+1. $$
\begin{thm}\label{dimthm}
The so defined algebra $\A(\D,\ga)$ is a Hopf algebra of dimension $N_1^{\binom{n_1+1}2}\cdots N_t^{\binom{n_t+1}2}.$
\end{thm}
\prf{The statement about the Hopf algebra is clear because we already know from earlier considerations in every component, that the ideal is a Hopf ideal. For the dimension we will use Theorem \ref{basisthm}. We have to check all the conditions.

We first need to prove the commutation relation between root vectors of one component and $\pth N$ powers of root vectors of the other components. And secondly, we have to show that the $u_{i,j}$ are central with regard to all generators $a_x$ of $\Uf(\D)$:
\begin{gather}
\label{commute}e_{S_k+i,S_k+j}^{\phantom{N_l}}e_{S_l+s,S_l+t}^{N_l}=\chi_{S_l+s,S_l+t}^{N_l}(g_{S_k+i,S_k+j})e_{S_l+s,S_l+t}^{N_l}e_{S_k+i,S_k+j}^{\phantom{N_l}}\\
\label{central}u_{S_k+i,S_k+j}a_x=a_xu_{S_k+i,S_k+j},
\end{gather}
where $1\le k,l\le t,\,1\le i<j\le n_k+1,\,1\le s<t\le n_l+1,\,1\le x\le S_t+n_t.$

When we are within one component the above equations follow from the original paper \cite{POIN}. 
The root vectors are linear combinations of homogeneous monomials in the generators $a_x$ of $\Uf(\D)$. Hence we see that Lemma \ref{techcom} will establish (\ref{commute}) for $k\neq l$.
 
(\ref{central}) is shown by induction on $j-i.$ Because of the recursive definition (\ref{defu}) of the $u_{i,j}$, the crucial part is:
\begin{align*}\ga_{k_{i,j}}(1-h_{S_k+i,S_k+j})a_x&=\ga_{k_{i,j}}a_x(1-\chi_x(h_{S_k+i,S_k+j})h_{S_k+i,S_k+j})\\&=\ga_{k_{i,j}}a_x(1-\chi_{S_k+i,S_k+j}^{N_k}(g_x^{-1})h_{S_k+i,S_k+j})\\&=a_x\ga_{k_{i,j}}(1-h_{S_k+i,S_k+j}).
\end{align*}
In the second step we used (\ref{cartype}) and as $x$ is not in the $\pth k$ component, all the corresponding entries of the Cartan matrix are zero. The third step uses the premise that $\ga$ is admissible.

Now we can apply Theorem \ref{basisthm} and the proof is finished.
}
\begin{lemma}\label{techcom} For all indices $i$ not in the $\pth k$ component of the diagram and $S_k<j\le l\le S_k+n_k+1$ we have with $q=\chi_j(g_j)=\chi_l(g_l)$\\

$i)\quad a_ie_{j,l+1}=$
\begin{subnumcases}{\hspace{-3em}=}
\chi_j(g_i)a_ja_i+\lam(1-g_ig_j)&\!if $j=l,$\label{case1}\\
\chi_{j,l+1}(g_i)e_{j,l+1}a_i+\lam(1-q^{-1})e_{j+1,l+1}&\!if $\lambda_{il}=0,\;j<l,$\label{case2}\\
\chi_{j,l+1}(g_i)e_{j,l+1}a_i-\lambda_{il}(1-q^{-1})\chi_{j,l}(g_i)e_{j,l}g_ig_l&\!otherwise.\label{case3}
\end{subnumcases}
$ii)\quad a_ie_{j,l+1}^{N_k}=\chi_{j,l+1}^{N_k}(g_i)e_{j,l+1}^{N_k}a_i.$
\end{lemma}

\prf{
$i)$ \bu The case $j=l$ is simply the defining relation (\ref{links}).

From now on $j<l$. We consider all possible linkings.

\bu If $i$ is not linked to any vertex $p$ with $i\le p\le l$, then $\lam=\lambda_{il}=0$ and a repeated use of (\ref{links}) gives (\ref{case2}).

\bu If $i$ is linked to $j$, then it can not be linked to $l$ as well. Hence we have to show the second case. We proceed by induction on $l-j$ and use the recursive definition (\ref{defroot}) of the root vectors.

For $l-j=1$ we have
\begin{align*} a_ie_{j,l+1}&=a_i[a_ja_l-\chi_l(g_j)a_la_j]\\&=\chi_j(g_i)a_ja_ia_l+\lam(1-g_ig_j)a_l-\chi_l(g_j)\chi_l(g_i)a_la_ia_j\\&=\chi_j(g_i)\chi_l(g_i)a_ja_la_i+\lam a_l-\lam\chi_l(g_ig_j)a_lg_ig_j\\&\quad-\chi_l(g_j)\chi_l(g_i)\chi_j(g_i)a_la_ja_i-\chi_l(g_j)\chi_l(g_i)a_l\lam(1-g_ig_j)\\&=\chi_{j,l+1}(g_i)[a_ja_l-\chi_l(g_j)a_la_j]a_i+\lam(1-\chi_l(g_i)\chi_l(g_j))a_l\\&=\chi_{j,l+1}(g_i)e_{j,l+1}a_i+\lam(1-\underbrace{\chi_i(g_l^{-1})\chi_j(g_l^{-1})}_{=1}\underbrace{\chi_j^{a_{jl}}(g_j)}_{=q^{-1}})a_l.
\end{align*}
We used (\ref{cartype}) and the condition $\chi_i\chi_j=1$ as $\lam\neq 0.$ For the induction step we use an analogue calculation. The last steps are as follows
\begin{align*} a_ie_{j,l+1}&=\chi_{j,l+1}(g_i)e_{j,l+1}a_i+\lam(1-q^{-1})[e_{j+1,l}a_l-\chi_l(g_{j,l})\chi_l(g_i)a_le_{j+1,l}]\\
&=\chi_{j,l+1}(g_i)e_{j,l+1}a_i+\\&\qquad+\lam(1-q^{-1})[e_{j+1,l}a_l-\chi_l(g_{j+1,l})\underbrace{\chi_l(g_j)\chi_l(g_i)}_{\chi_j(g_l^{-1})\chi_i(g_l^{-1})}a_le_{j+1,l}]\\&=\chi_{j,l+1}(g_i)e_{j,l+1}a_i+\lam(1-q^{-1})e_{j+1,l+1}.
\end{align*}
\bu If $i$ is linked to $l$ we have $\lambda_{i,l}\neq 0$ and hence we need to prove (\ref{case3}). A direct calculation gives
\begin{align*}
a_ie_{j,l+1}&=a_i[e_{j,l}a_l-\chi_l(g_{i,l})a_le_{j,l}]\\
&=\chi_{j,l}(g_i)e_{j,l}a_ia_l-\chi_l(g_{j,l})\chi_l(g_i)a_la_ie_{j,l}-\chi_l(g_{j,l})\lambda_{il}(1-g_ig_l)e_{j,l}\\
&=\chi_{j,l}(g_i)e_{j,l}\chi_l(g_i)a_la_i+\chi_{j,l}(g_i)e_{j,l}\lambda_{il}(1-g_ig_l)\\
&\quad -\chi_l(g_{j,l})\chi_l(g_i)a_l\chi_{j,l}(g_i)e_{j,l}a_i-\chi_l(g_{j,l})\lambda_{il}e_{j,l}\\
&\quad +\chi_l(g_{j,l})\lambda_{il}\chi_{j,l}(g_ig_l)e_{j,l}g_ig_l\\
&=\chi_{j,l+1}(g_i)[e_{j,l}a_l-\chi_l(g_{i,l})a_le_{j,l}]a_i+[\chi_{j,l}(g_i)-\chi_l(g_{j,l})]\lambda_{il}e_{j,l}\\
&\quad +[\chi_l(g_{j,l})\chi_{j,l}(g_ig_l)-\chi_{j,l}(g_i)]\lambda_{il}e_{j,l}g_ig_l.
\end{align*}
As $i$ is not in the $\pth k$ component we have by (\ref{cartype}) and $\chi_i\chi_l=1$
$$\chi_{j,l}(g_i)=\chi_i^{-1}(g_{j,l})=\chi_l(g_{j,l}).$$
Hence the second term in the last step of the above calculation vanishes, and for the bracket of the third term we calculate
\begin{align*}
\chi_l(g_{j,l})\chi_{j,l}(g_ig_l)-\chi_{j,l}(g_i)&=\chi_{j,l}(g_i)(\chi_l(g_{j,l})\chi_{j,l}(g_l)-1)\\
\chi_l(g_{j,l})\chi_{j,l}(g_l)&=\chi_l(g_j)\chi_l(g_{j+1})\cdots\chi_l(g_{l-1})\\
&\quad\cdot\chi_j(g_l)\chi_{j+1}(g_l)\cdots\chi_{l-1}(g_l)\\
&=1\cdot1\cdots\chi_l(g_l)^{-1}.
\end{align*}
\bu For the last case where $i$ is linked to a vertex $p$ with $j<p<l,$ we again proceed by induction on $l-p.$
As $\lam=\lambda_{jl}=0$ we have to show (\ref{case2}).

If $l-p=1$ we use the recursive definition of the root vectors and then (\ref{case3}). We set $F:=\lambda_{il-1}(1-q^{-1})\chi_{j,l-1}(g_i)$ and have
\begin{align*}
a_ie_{j,l+1}&=a_i[e_{j,l}a_l-\chi_l(g_{i,l})a_le_{j,l}]\\
&=[\chi_{j,l}(g_i)e_{j,l}a_i+F\cdot e_{j,l-1}g_ig_{l-1}]a_l\\
&\quad -\chi_l(g_{i,l})\chi_l(g_i)a_l[\chi_{j,l}(g_i)e_{j,l}a_i+F\cdot e_{j,l-1}g_ig_{l-1}]\\
&=\chi_{j,l}(g_i)e_{j,l}\chi_l(g_i)a_la_i-\chi_l(g_{i,l})\chi_l(g_i)\chi_{j,l}(g_i)a_le_{j,l}a_i\\
&\quad +F\cdot\chi_l(g_ig_{l-1})e_{j,l-1}a_lg_ig_{l-1}-F\cdot\chi_l(g_{i,l})\chi_l(g_i)a_le_{j,l-1}g_ig_{l-1}\\
&=\chi_{j,l+1}(g_i)e_{j,l+1}a_i\\&\quad+F\cdot\chi_l(g_ig_{l-1})[e_{j,l-1}a_l-\chi_l(g_{i,l-1})a_le_{j,l-1}]g_ig_{l-1}.
\end{align*}
However, the last square bracket is zero according to \cite[(7.17)]{POIN}. The induction step is now simple.\\

$ii)$ \bu The case $\lam=\lambda_{il}=0$ is trivial.

\bu So let now $j=l$ and $\lam\neq 0$. Then $\chi_j(g_i)=\chi_i^{-1}(g_j)=\chi_j(g_j)=q$ and using (\ref{case1}) we have
\begin{align*}a_ia_j^{N_k}&=q^{N_k}a_j^{N_k}a_i+\lam(1+q+q^2+\dots+q^{N_k-1})a_j^{N_k-1}(1-q^{N_k-1}g_ig_j)\\
&=a_j^{N_k}a_i.\end{align*}
Here we used the fact that $N_k$ is the order of $q$.\\
From now on again $j<l.$

\bu If $\lam\neq 0,$ we set $x=e_{j,l+1},$ $y=a_i,$ $z=\lam(1-q^{-1})e_{j+1,l+1},$ $\alp=\chi_{j,l+1}(g_i)$ and $\beta=\chi_{j+1,l+1}^{-1}(g_{j,l+1}).$ Then, because of \cite[(7.24)]{POIN} $zx=\beta xz.$ Hence, cf. \cite[Lemma 3.4 (1)]{ASa2},
$$yx^{N_k}=\alp^{N_k}x^{N_k}y+\left(\sum_{m=0}^{N_k-1}\alp^m\beta^{N_k-1-m}\right)x^{N_k-1}z.$$
Using $\chi_i\chi_j=1$ and (\ref{cartype}) we see that $\alp=\chi_i^{-1}(g_{j,l+1})=\chi_j(g_{j,l+1})$ and so
\begin{align*}
\alp^m\beta^{N_k-1-m}&=\beta^{N_k-1}\chi_j^m(g_{j,l+1})\chi_{j+1,l+1}^m(g_{j,l+1})\\
&=\beta^{N_k-1}\chi_{j,l+1}^m(g_{j,l+1})=\beta^{N_k-1}(B_{j,l+1}^{j,l+1})^m=\beta^{N_k-1}q^m.
\end{align*}
The last equality follows from \cite[(7.5)]{POIN}. The geometric sum gives zero again.

\bu The final case $\lambda_{i,l}\neq 0$ is treated similarly to the previous one. This time $z=-\lambda_{il}(1-q^{-1})\chi_{j,l}(g_i)e_{j,l}g_ig_l$ and $\beta=\chi_{j,l+1}(g_ig_l)\chi_{j,l+1}(g_{j,l}),$ because of \cite[(7.23)]{POIN}. So we have
\begin{align*}
\alp^m\beta^{N_k-1-m}&=\beta^{N_k-1}\chi_{j,l+1}^m(g_i)\chi_{j,l+1}^{-m}(g_i)\chi_{j,l+1}^{-m}(g_{j,l+1})\\
&=\beta^{N_k-1}(B_{j,l+1}^{j,l+1})^{-m}=\beta^{N_k-1}q^{-m}.
\end{align*}
}
Here is the final result.
\begin{thm}
Let $\D$ and $\D'$ be two linking data as above and $\ga$ and $\ga'$ two admissible parameter families. Then $\A(\D,\ga)$ and $\A(\D',\ga')$ are quasi-isomorphic.
\end{thm}
The proof is just a combination of the results obtained in the previous sections.

First, one shows that the Hopf algebras are quasi-isomorphic to ones where all the parameters $\ga$ are zero. This is done component-wise and each step is identical to the proof of Theorem \ref{main}. The necessary commutation relations are in Lemma \ref{techcom}.

This way we find that the Hopf algebras $\A(\D,\ga)$ and $\A(\D',\ga')$ are quasi-isomorphic to $\uf(\D)$ and $\uf(\D'),$ respectively. Now Theorem \ref{linkthm} completes the proof.\\

It should be possible to extend this approach to arbitrary Dynkin diagrams of finite Cartan type. The only crucial part are the commutation relations between the root vectors. They should be checked for the other diagrams.


\begin{thebibliography}{abcd}

\bibitem[AS1]{P^3} \liti{N.~Andruskiewitsch, H.-J.~Schneider}{Lifting of Quantum Linear Spaces and Pointed Hopf algebras of Order $p^3$}{J.~Algebra}{209}{1998}{658-691}
\bibitem[AS2]{ASa2} \liti{N.~Andruskiewitsch, H.-J.~Schneider}{Lifting of Nichols Algebras of Type $A_2$ and Pointed Hopf Algebras of order $p^4$}{ in ``Hopf algebras and quantum groups'', Proceedings of the Brussels Conference, eds. S. Caenepeel, F. Van Oystaeyen, Lecture Notes in Pure and Appl. Math.}{209}{2000}{1-14, Marcel Dekker, New York}
\bibitem[AS3]{AS} \liti{N.~Andruskiewitsch, H.-J.~Schneider}{Finite quantum groups over abelian groups of prime exponent}{Ann. Sci. Ec. Norm. Super.}{35}{2002}{1-26}
\bibitem[AS4]{POIN} \liti{N.~Andruskiewitsch, H.-J.~Schneider}{Pointed Hopf algebras}{}{}{2002}{to appear in "New directions in Hopf algebras", MSRI series Cambridge Univ. Press}
\bibitem[AS5]{New-AS} \liti{N.~Andruskiewitsch, H.-J.~Schneider}{A characterization of Quantum Groups}{}{}{2002}{available as {\tt math.QA/0201095}}
\bibitem[BDR]{BDR} \liti{M.~Beattie, S.~D\u asc\u alescu, S.~Raianu}{Lifting of Nichols algebras of type $B_2$}{preprint}{}{2001}{available at\\ {\tt http://www.mta.ca/\~{}mbeattie/research/prepr.htm}}
\bibitem[D]{D} \liti{D.~Didt}{Linkable Dynkin diagrams}{preprint}{}{2001}{available as {\tt math.QA/0201268}}
\bibitem[EG]{EG} \liti{P.~Etingof, S.~Gelaki}{On Families of Triangular Hopf algebras}{preprint}{}{2001}{available as {\tt math.QA/0110043}}
\bibitem[M]{MA} \liti{A.~Masuoka}{Defending the negated Kaplansky conjecture}{Proc. Amer. Math. Soc.}{129}{2001}{3185-3192}
\bibitem[S]{S} \liti{P.~Schauenburg}{Hopf bi-Galois extensions}{Comm. Algebra}{24}{1996}{3797-3825}


\end{thebibliography}
\end{document}